\documentclass[12pt]{amsart}
\usepackage{amsmath,amsthm,amscd,amsfonts,amssymb}
\usepackage[all]{xy}

\begin{document}
   \baselineskip   .7cm
\newtheorem{defn}{Definition}[section]
\newtheorem{thm}{Theorem}
\newtheorem{propo}[defn]{Proposition}
\newtheorem{cor}[defn]{Corollary}
\newtheorem{lem}[defn]{Lemma}
\theoremstyle{remark}
\newtheorem{rem}{Remark}[section]
 \newtheorem{ex}[defn]{Example}
\renewcommand\o{{{\mathcal O}}}
\newcommand\s{\sigma}
\newcommand\w{\widehat}
\newcommand\cal{\mathcal}

\title [The p-rank of the reduction $\rm{mod}\, p$   of Jacobians]
{The $p$-rank of the reduction $\rm{mod}\, p$     of Jacobians and Jacobi sums}

\author[A. \'Alvarez]{A. \'Alvarez${}^*$     }
\address{Departamento de Matem\'aticas  \\
Universidad de Salamanca \\ Plaza de la Merced 1-4. Salamanca
(37008). Spain.}

\thanks{MSC:   11G30, 11G20, 11G10 \\$*$  Departamento de Matem\'aticas.
Universidad de Salamanca. Spain.}

\maketitle

  \begin{abstract} Let $Y_K\to X_K$ be a ramified cyclic covering of curves,  where $K$ is a cyclotomic field. In this work we study the $p$-rank of the reduction ${\rm{mod}\, p}$  of a model of the Jacobian of $Y_K$.  In this way, we obtain counterparts   of the Deuring polynomial, defined for elliptic curves, for genus greater than one.   To carry out this study we use the relationship between Jacobi sums and   $L$-functions. This is  established  in \cite{W} for the case of Fermat curves.

  \end{abstract}

\tableofcontents


\section{Introduction and previous notation}

{\bf{Previous notation:}}
Let $m$ be a prime integer and $K:={\Bbb Q}(\epsilon_m)$ the cyclotomic field with $\epsilon_m$ an $m$-primitive root  of the unity, and let   $t$ be an integer prime to $m$. We  consider the automorphism $\sigma_t$ of $K$ defined by $\sigma_t(\epsilon_m):=\epsilon_m^t$. Let  $X_K$   be a proper and geometrically irreducible curve over $K$, of genus $g$. Let ${  x}_0,\cdots,{  x}_d$ be different points within $X_K$.

Now  let us  consider $Y_K\to X_K$, a Galois ramified covering,  of Galois group $G:={\Bbb Z}/m$,   ramified at
${  x}_0,\cdots,{  x}_d$. Let $\Sigma_{Y_K}$, $\Sigma_{X_K}$ be the function fields of $Y_K$ and $X_K$, respectively. We have that $\Sigma_{Y_K}=\Sigma_{X_K}(\sqrt[m]{ f })$, where $f\in \Sigma_{X_K}$. We assume that $f$ is not a $m$-power and thus that $ {Y_K}$ is connected. Let ${\rm{div}}(f )$ be the principal divisor associated with $f \in \Sigma_{X_K}$. By Kummer theory we can choose $f$ such that
${\rm{div}}(f)=a_0\cdot x_0+\cdots+a_d\cdot x_d+m\cdot D$, with $0<  a_i  <m$ and $D$ is a divisor on $X_K$. Note that $a_0  +\cdots+a_d =0\,{\rm{mod}}\, m$.

We choose $n$, a positive integer, such that $X_K$ has good reduction, $X$, over ${\rm{Spec}}(A)$ with  $A:={\Bbb Z}[\epsilon_m,\frac{1}{m\cdot n}]$  and that the points ${  x}_0,\cdots,{  x}_d$ do not coalesce $\text{ mod } \mathfrak p$, for each $  \mathfrak p \in {\rm{Spec}}(A)$. With these conditions there exists   a model  $Y\to X$ over ${\rm{Spec}} (A)$ for $Y_K\to X_K$, c.f.\cite{Be} .

We denote by $X_{\mathfrak p}$  the reduction at ${\mathfrak p}\in {\rm{Spec}} (A)$ of $X$ and by $Y_{\mathfrak p}$ a proper, smooth model for the reduction at ${\mathfrak p}\in {\rm{Spec}} (A)$  of    $Y$. We denote by $\Sigma_{Y_\mathfrak p}$, $\Sigma_{X_\mathfrak p}$ the function fields of $Y_\mathfrak p$ and $X_\mathfrak p$, respectively. Let us denote  $k({\mathfrak p})$  the residual field of ${\mathfrak p}$, which we assume to have $ q=p^h$ elements.
We also assume that $Y_{\mathfrak p}\to X_{\mathfrak p}$ is a ramified Galois covering of group $G$, ramified at
$\overline x_0,\cdots,\overline x_d$, the reduction of $x_0,\cdots,x_d$ at $\mathfrak p$ and that
  $\Sigma_{Y_\mathfrak p}=\Sigma_{X_\mathfrak p}(\sqrt[m]{    f})$ where ${\rm{div}}(f)=a_0\cdot \overline x_0+\cdots+a_d\cdot \overline  x_d+m\cdot \overline  D$, with   $\overline  D$   a divisor on $X_\mathfrak p$.
We denote $T:=\{\overline x_0,\cdots, \overline  x_d\}$.

We  fix  an algebraic closure ${\Bbb F}$ of ${\Bbb F}_p$. If $Z$ is a variety over ${\Bbb F}_p$ we denote $\overline Z:=Z\otimes_{{\Bbb F}_p}{\Bbb F}$.
 In this article,
we study the  characteristic polynomials of the $p^h$-Frobenius morphism, $F_\mathfrak p$, of certain  ${\Bbb Q}_l[\epsilon_m]$-modules associated with the $l$-adic cohomology group
  $ H^1(\overline {Y}_{\mathfrak p} ,{\Bbb Q}_l) $.  These polynomials are     $L$-functions   of $X_{\mathfrak p}$.

The first part of this article is devoted to proving that the constant term of these polynomials is given by Jacobi sums. This   has already been  proved in \cite{W} for $g=0$, $d=2$ and in general    in \cite{D}, bearing in mind that these terms are essentially the  constant functions of certain Dirichlet $L$-series. In this work,   we use geometric  methods to achieve these results.

The second part of this article addresses   the $p$-rank of ${\rm{Pic}}^0_{Y_{\mathfrak p}}$.
We say that an abelian variety, $\mathcal A$, defined  over ${\Bbb F}_\mathfrak p$, is supersingular when $\overline {\mathcal A} $ is isogenous to a product of supersingular elliptic curves. Equivalently,  $\mathcal A$ is supersingular when the eigenvalues of the Frobenius morphism  $F_{\mathfrak p}$ are $\zeta \cdot q^{1/2}$, $\zeta $ being a root of the unity.

Let us denote by $[\frac{a}{b}]$   the integer part of the fraction $\frac{a}{b}$ and  $<\frac{a}{b}>:=\frac{a}{b}-[\frac{a}{b}]$.
By using   properties of Jacobi sums we prove that if   $p$ splits completely in  ${\Bbb Z}[\epsilon_m, \frac{1}{m\cdot n}]$ and there exists      $t\in \{ 1,\cdots, m-1\}$  with
$$[<\frac{a_1}{m}>+\cdots +<\frac{a_d}{m}>]\neq [<\frac{t\cdot a_1}{m}>+\cdots +<\frac{t\cdot a_d}{m}>].$$
Then  the Jacobian of $Y_{\mathfrak p}$ is not a supersingular abelian variety.

Let $[p]$ be  the multiplication by $p$ on ${\rm{Pic}}^0_{Y_{\mathfrak p}}$. We say that $Y_{\mathfrak p}$ has $p$-rank $0$  when   $ {\rm{Ker}} [p]_{{\rm{red}}}={\rm{Spec}}(k(\mathfrak p))$, or equivalently $ {\rm{Ker}} [p]({\Bbb F})=\{0\}$.
We   prove  for $X_{\mathfrak p}={\Bbb P}^1$ that
 if   $t\in \{ 1,\cdots, m-1\}$  with
$$  [<\frac{t\cdot a_1}{m}>+\cdots +<\frac{t\cdot a_d}{m}>]= 0,$$
     then $Y_{\mathfrak p}$ does not have $p$-rank $0$.

Let  $p$  be inert in  ${\Bbb Z}[\epsilon_m, \frac{1}{m\cdot n}]$ and let us  consider   the integers  $a_1,\cdots, a_d$   satisfying  $0<  a_i  <m$ for each $1\leq i \leq d-2$, and $a_1+\cdots+a_{d }\neq 0$ mod $m$. We prove that    the proper, smooth model of the curve
$$y^m-x^{a_1}(x-1)^{a_2}(x-\alpha_1)^{a_3}\cdots (x-\alpha_{d-2})^{a_{d }},$$
 defined on $k(\mathfrak p)={\Bbb F}_{p^{m-1}},$ has $p$-rank $0$ if and only if $\alpha_1,\cdots, \alpha_{d-2}$ satisfy
a system of $d-2$ algebraic  equations    of degrees $\leq l\cdot(a_3+\cdots+a_d)(\frac{q^{m-1}-1}{m})$ with $1\leq l \leq d-2$, respectively. In \cite[5]{Bo} the generalized Hasse-witt matrix it is obtained explicitly, and hence the   equations to study the $p$-rank of ${\rm{Pic}}^0_{Y_{\mathfrak p}}$.

 When $g=0$, $d=3$, $a_1=a_2=a_3=1$ and $m=2$,  one obtains the   Deuring polynomial  whose roots are  the values such that the elliptic curve $y^2-x(x-1)(x-\lambda)$ is supersingular,
$$H(\lambda)= \sum^n_{i=0}\left(\begin{matrix}    r
 \\ i  \end{matrix}\right)^2\lambda^i,(\text{with} \quad r=(p-1)/2).$$
 We also study when $Y_{\mathfrak p}$ has $p$-rank $0$  in the case of $p$ splitting completely in ${\Bbb Z}[\epsilon_m, \frac{1}{m\cdot n}]$.

In \cite{Bo}  the Hasse-Witt invariants of $Y_\mathfrak p$ are studied and calculated. These invariants give an upper bound for the $p$-rank of $Y_\mathfrak p$. It is showed that for $p$ large  and  $\overline x_0,\cdots,\overline x_d$ generic points, this upper bound is equal to the $p$-rank. In \cite{E},
for $X_\mathfrak p={\Bbb P}^1$,  the author gives bounds for the $p$-rank of $Y_\mathfrak p$. In \cite{G}  the Hasse-Witt matrix is calculated to show that for the Fermat curves there exists a set of primes $p$ with positive density such that the Fermat curves $(\text {mod } p)$ are not supersingular but their $p$-rank is $0$. The dimension of the moduli of hyperelliptic curves with fixed $p$-rank is obtained in \cite{GlP}. For Fermat curves in \cite {N}  the supersingularity of these curves is studied. In \cite{Y} The $p$-rank of a  hyperelliptic curve is given  by the  rank of the Cartier-Manin matrix.

\section{Cyclic extensions}

We now give some general notation that we shall use along this work.
Let $s$ be a global section of a line bundle $L$ on $X_\mathfrak p$. We also  denote by $s$ the morphism, of  ${\mathcal O}_{X_\mathfrak p}$-modules, ${\mathcal O}_{X_\mathfrak p} \to L$,  such that $1\to s$. Here, ${\mathcal O}_{X_\mathfrak p}$ denotes the   sheaf of rings associated with $X_\mathfrak p$.

We denote by $E$ and ${\mathfrak m}$  the effective divisor $\overline x_0+\cdots+\overline  x_d$ on $X_\mathfrak p$  and the ideal associated with $  E$ inside ${\mathcal O}_{X_\mathfrak p}$, respectively.

If $z\in k(\mathfrak p)^\times={\Bbb F}^\times_{p^h}$, then we denote by $\chi_{\mathfrak p}(z)$ the unique $m$-root of the unity such that $\chi_{\mathfrak p}(z)=z^{\frac{p^h-1}{m}}\, {\rm{mod}}\, \mathfrak p$.

Let $L$ be a line bundle over $X_\mathfrak p$ and $\iota_\mathfrak m:L\to {\mathcal O}_{X_\mathfrak p}/{\mathfrak m}$ a surjective morphism of ${\mathcal O}_{X_\mathfrak p}$-modules.
Let ${\rm{Pic}}^0_{X_\mathfrak p,{\mathfrak m} }$ be the generalized Jacobian for   $\mathfrak m$. We have that ${\rm{Pic}}^0_{X_\mathfrak p,{\mathfrak m} }$ is a scheme over $k(\mathfrak p)$ that represents  isomorphism classes of pairs $(L,\iota_\mathfrak m)$ ($\mathfrak m$-level structures). We say that two level structures $(L,\iota_\mathfrak m)$ and $(L',\iota'_\mathfrak m)$ are equivalent  when there exists an isomorphism of line bundles $ u:L\to L'$ such that $\iota'_\mathfrak m\cdot u =\iota_\mathfrak m$.

Let us consider $\Sigma^{\mathfrak m}_{X_\mathfrak p}:=\{g\in \Sigma_{X_\mathfrak p}\text{ such that\,} g\equiv 1\, \text{mod }\, {\mathfrak m}\}$. The  equivalence classes of level structures  are in  one-to-one correspondence with the $\mathfrak m$-equivalence classes of divisors on $X_\mathfrak p$ supported outside $T$;  two divisors $D $ and $D'$  are ${\mathfrak m}$-equivalent when there exists a $g\in \Sigma^{\mathfrak m}_{X_\mathfrak p}$ with  $D -D'={\rm{div}}(g)$.

Let $\pi$ be the natural epimorphism ${\mathcal O}_{X_\mathfrak p}\to {\mathcal O}_{X_\mathfrak p}/{\mathfrak m}$. We shall use the term    (projective) space of $\mathfrak m$-sections of a level structure $(L,\iota_\mathfrak m)$, to the subspace $H^0_\mathfrak m((L,\iota_\mathfrak m))\subset H^0(X_\mathfrak p,L ) $,    of global sections of $L$, $s:{\mathcal O}_{X_\mathfrak p}\to L$ such that $\iota_{\mathfrak m}\cdot s=\pi$. The effective  $\mathfrak m$-equivalent divisors associated with $(L,\iota_\mathfrak m)$ are given by the zero locus of the $\mathfrak m$-sections   $s$; Let us consider  inclusions $L, L' \subset \Sigma_{X_\mathfrak p}$. We have that   $(L',\iota'_\mathfrak m)$ is equivalent to $(L,\iota_\mathfrak m))$ if and only if there exists a $g\in {\Sigma}^{\mathfrak m}_{X_\mathfrak p}$ verifying
$g\cdot L'=L$.
 Because the difference between two $\mathfrak m$-sections of $(L,\iota_\mathfrak m)$ is a global section of   $L(-E)$, we have that if $s\in H^0_\mathfrak m((L,\iota_\mathfrak m)) $, then $H^0_\mathfrak m((L,\iota_\mathfrak m))=s+H^0(X_\mathfrak p,L(-E))$, .

We   denote  by   $I_T$, $I^0_T$ and $O_T^\times$ the ideles,  ideles of degree $0$ and integer ideles on $\Sigma_{X_\mathfrak p}$ outside $T$, respectively.

We consider the $p^h$-Frobenius morphism, $F_\mathfrak p$,
 and the Lang isogeny  $P:=F_\mathfrak p-Id:{\rm{Pic}}^0_{X_\mathfrak p,{\mathfrak m} }\to {\rm{Pic}}^0_{X_\mathfrak p,{\mathfrak m} }$. Bearing in mind a divisor of degree $1$ supported on $\vert X_\mathfrak p\vert \setminus T$,  $D_1$, we have, by translation, an immersion $X_{\mathfrak p}\setminus T \to {\rm{Pic}}^0_{X_\mathfrak p,{\mathfrak m} }$.
 By class field theory for function fields over finite fields, \cite[VI, nº11, Proposition 9; VI, nº6, Proposition 6 and VI, nº1, Proposition 11]{S}, we have that
  $P^{-1}(X_{\mathfrak p}\setminus T)$ gives  the ${\mathfrak m}$-ray class field $  H_\mathfrak m$ for $\Sigma_{X_\mathfrak p}$. It is a Galois extension with Galois group isomorphic to ${\rm{Pic}}^0_{X_\mathfrak p,{\mathfrak m} }(k(\mathfrak p)) =\frac{ I^0_T}{(\Sigma^{\mathfrak m}_{X_\mathfrak p})^\times\cdot O_T^\times}$ . This   group is the group of the ${\mathfrak m}$-equivalence classes of divisors on $X_{\mathfrak p}$ supported outside $T$.

Therefore, the cyclic extension  $\Sigma_{Y_\mathfrak p}/\Sigma_{X_\mathfrak p}$ of Galois group ${\Bbb Z}/m$,  where  $\Sigma_{Y_\mathfrak p}=\Sigma_{X_\mathfrak p}(\sqrt[m]{    f})$, gives an   epimorphism of groups
$$ \frac{ I^0_T}{(\Sigma^{\mathfrak m}_{X_\mathfrak p})^\times\cdot O_T^\times}\to {\Bbb Z}/m.$$
This morphism is given by   the Artin map \cite[VI, nº29, Th\'eor\`eme 9]{S}
$$(\quad ,\Sigma_{Y_\mathfrak p}/\Sigma_{X_\mathfrak p} ):\frac{ I_T}{(\Sigma^{\mathfrak m}_{X_\mathfrak p})^\times\cdot O_T^\times}\to {\Bbb Z}/m.$$
 If $\nu_1$ is the idele class  of degree $1$ associated with the  divisor $D_1$, then
 $(\nu_1,\Sigma_{Y_\mathfrak p}/\Sigma_{X_\mathfrak p})=1.$

By noting that the residual fields $k(\overline x_i)$, ($i=0,\cdots,d$),  are isomorphic to $k(\mathfrak p)$, the morphism of forgetting the level structure, $(L,\iota_\mathfrak m)\to L$, gives the exact sequence of schemes in groups
$$1\to  {\Bbb G}_m  \times \overset {\underset \smile {d+1}}\cdots \times {\Bbb G}_m/{\Bbb G}_m \to {\rm{Pic}}^0_{X_\mathfrak p,{\mathfrak m} }\to  {\rm{Pic}}^0_{X_\mathfrak p} \to 1$$
and thus we have an injective morphism of groups,
$$  \eta:k(\overline x_0)^\times \times\cdots\times k(\overline x_d)^\times/k({\mathfrak p})^\times  \to \frac{ I^0_T}{(\Sigma^{\mathfrak m}_{X_\mathfrak p})^\times\cdot O_T^\times} .$$

If we consider $\sigma \in \frac{I^0_T}{(\Sigma^{\mathfrak m}_{X_\mathfrak p})^\times\cdot O_T^\times}$,  then via the Artin symbol $\sigma(\sqrt[m]{    f})=z\cdot (\sqrt[m]{    f})  $ , with $z\in k(\mathfrak p)^\times     ={\Bbb F}^\times_{p^h}$.
We denote by $\chi_f$ the character of $G$ defined by $\chi_f(\sigma):=\chi_{\mathfrak p}(z)$.

For $(z_0,\cdots,z_d)\in k(\overline x_0)^\times \times\cdots\times k(\overline x_d)^\times/k({\mathfrak p})^\times$ we have that $$\chi_f( \eta(z_0,\cdots,z_d) ,\Sigma_{Y_\mathfrak p}/\Sigma_{X_\mathfrak p})=\chi^{-a_0}_{\mathfrak p}(z_0)\cdots \chi^{-a_d}_{\mathfrak p}(z_d).$$

Note that  we have an isomorphism $$k(\overline x_0)^\times \times\cdots\times k(\overline x_d)^\times/k({\mathfrak p})^\times\simeq k(\overline x_1)^\times\times\cdots\times k(\overline x_d)^\times,  $$
  and therefore we   can assume that $z_0=1$.

\section{L-functions}
This section is devoted to studying certain incomplete $L$-functions of the curve $X_{\mathfrak p}$ over $k({\mathfrak p})$, we follow \cite{A} and \cite{T}. The action of $ G$ over $ H^1(\overline {Y}_{\mathfrak p} ,{\Bbb Q}_l) $ gives a decomposition into eigenspaces
$$H^1(\overline {Y}_{\mathfrak p} ,{\Bbb Q}_l) =\oplus_{\chi}H^1(\overline {Y}_{\mathfrak p} ,{\Bbb Q}_l)^{\chi}.$$
The sum is over all characters of $ G:={\Bbb Z}/m $.
In this section  we shall   calculate the characteristic polynomial of $F_\mathfrak p$ as an endomorphism of the ${\Bbb Q}_l[\epsilon_m]$-module $H^1(\overline {Y}_{\mathfrak p} ,{\Bbb Q}_l)^{\chi_f}$.

Let $F_x$ be the Frobenius element, with $  x\in \vert X_\mathfrak p\vert \setminus T$. We consider it      as   an element of    the Galois group,  $G_{\mathfrak m}:=\frac{ I^0_T}{(\Sigma^{\mathfrak m}_{X_\mathfrak p})^\times\cdot O_T^\times}$, of the ray class field   $H_\mathfrak m/\Sigma_{X_\mathfrak p}$. Let $t_x$ be a local parameter for $x$.
We have that $F_x=(t_x, H_\mathfrak m/\Sigma_{X_\mathfrak p} ) $, where $(\quad , H_\mathfrak m/\Sigma_{X_\mathfrak p} )$ is the Artin symbol for the Galois extension $H_\mathfrak m/\Sigma_{X_\mathfrak p}$.
We now consider the $T$-incomplete $L$-function
$$\theta_{H_\mathfrak m/\Sigma_{X_\mathfrak p}, T}(t):=\underset{x\in \vert X_\mathfrak p\vert \setminus T}\prod (1-F_x\cdot t^{deg(x)})^{-1}.$$

Let  $N$  be  a divisor supported outside $T$ and  with class $[N]\in G_{\mathfrak m}$.
We denote by $\sigma_N$ the element of the Galois group of the extension $ H_\mathfrak m/\Sigma_{X_\mathfrak p}$ associated with $N$  via the Artin symbol.

Let denote us by   $L_{\mathfrak m}(N,i )$   the cardinal of the set of effective divisors on $X_\mathfrak p$ supported
outside
$T$   and
${\mathfrak m}$-equivalents to $N+i\cdot D_1 $.

Let $D=n_1\cdot y_1+\cdots+n_r\cdot  y_r$ be a divisor on $X_\mathfrak p $ with support outside $T$,  with $t_{y_i}$    local parameters for $y_i$, $(i=1,\cdots,r)$. We define
$$(D , H_\mathfrak m/\Sigma_{X_\mathfrak p} ):=\prod_{i=1}^r  (t_{y_i},H_\mathfrak m/\Sigma_{X_\mathfrak p})^{n_i}$$

Similar to \cite{A} 4.1.1, we can compute this $L$-function in terms of $G_{\mathfrak m}$:
$$\theta_{H_\mathfrak m/\Sigma_{X_\mathfrak p},T}(t)
 =\sum_{[N]\in G_{\mathfrak m}}  \sigma_N\cdot (  \underset{i=0}{\overset
{ 2g+d-1} \sum} L_{\mathfrak m}(N,i )\cdot t^{i}
 + \underset{j\geq 0}  \sum q^{g+j}\cdot t^{ 2g+d+j}) =    $$
 $$=  \underset{i=0}{\overset
{ 2g+d-1} \sum} \sum_{[N]\in G_{\mathfrak m}} (L_{\mathfrak m}(N,i )\cdot  \sigma_N)\cdot t^{i}
 +  ( \sum_{[N]\in G_{\mathfrak m}}  \sigma_N)(\frac{q^g\cdot t^{2g+d}}{1-qt})=$$
 $$=\underset{i=0}{\overset
{ 2g+d-1} \sum}\sum_{D, {\rm{deg}}(D)=i}     (D , H_\mathfrak m/\Sigma_{X_\mathfrak p} ) \cdot t^{i}+  ( \sum_{[N]\in G_{\mathfrak m}}  \sigma_N)(\frac{q^g\cdot t^{2g+d}}{1-qt}), $$
where the sum is over all effective divisors $D$ on $X_\mathfrak p$ with support outside $T$.

 Note that if $ (L,\iota_\mathfrak m) $ is  a  level structure associated with $N+i\cdot D_1 $ then $L_{\mathfrak m}(N,i )=\#H^0_\mathfrak m((L,\iota_\mathfrak m)) $.  Thus, either $L_{\mathfrak m}(N,i )=0$ or $L_{\mathfrak m}(N,i )=\#H^0(X_{\mathfrak p}, N+i\cdot D_1 -E)$. The divisor $D_1$ is defined in section $2$.

Let $y$ be a geometric point of $X_\mathfrak p$. We denote by $f(y)$ the value of $f$ within the residual field $k(y)\subset {\Bbb F}$.

Let $D=n_1\cdot y_1+\cdots+n_r\cdot y_r$ be a divisor on $X_\mathfrak p $, we define
$$f(D):=\prod_{i=1}^r f(y_i)^{n_i \cdot \frac{p^{h\cdot {\rm{deg}}(y_i)-1}}{p^h-1}}  ,$$
 If $x\in \rm{sup}({\rm{div}}(f))\setminus T$ then we define $f(x)$  by considering a divisor $D'$ linearly $\mathfrak m$-equivalent to $x$ and ${\rm{sup}}(D')\cap {{\rm{sup}}}({\rm{div}}(f))=\emptyset$.

Let $P^{\chi_f}_{F_\mathfrak p}(t):={\rm{det}} (t-{F_\mathfrak p} )\in {\Bbb Z}[\epsilon_m][t] $ be the characteristic polynomial of the endomorphism $F_\mathfrak p$ over the ${\Bbb Q}_l(\epsilon_m)$-vector space $H^1(\overline {Y}_{\mathfrak p} ,{\Bbb Q}_l)^{\chi_f}$.

 By \cite[3.5]{T}, we have that this polynomial is equal to
$$ (*)  \chi_f(t^{2g+d-1}\cdot \theta_{H_\mathfrak m/\Sigma_{X_\mathfrak p},T} (\frac{1}{t}  ))=   \underset{i=0}{\overset
{ 2g+d-1} \sum} \sum_{[N]\in G_{\mathfrak m}}   L_{\mathfrak m}(N,i )\cdot\chi_f(\sigma_N)\cdot t^{2g+d-1-i}  $$
$$=  \underset{i=0}{\overset
{ 2g+d-1} \sum}\sum_{D, {\rm{deg}(D)}=i}    \chi_{\mathfrak p}(f(D))t^{2g+d-1-i}, $$
where $\chi_f $ is defined in section 2 and the sum is over all effective divisors $D$   on $X_\mathfrak p$ with support outside $T$.

If $p_0(t)$ is the characteristic polynomial of $F_\mathfrak p$ as an endomorphism of the ${\Bbb Q}_l$-vector space
 $H^1(\overline X_{\mathfrak p},{\Bbb Q}_l) $, then
 the characteristic polynomial of $F_\mathfrak p$, where $F_\mathfrak p$ is considered as an endomorphism of the ${\Bbb Q}_l $-vector space $H^1(\overline Y_{\mathfrak p},{\Bbb Q}_l) $, is
$$ p_0(t)\cdot \underset {1\leq j <m }\prod P^{\chi^j_f}_F(t)\in {\Bbb Z}[t].$$

\bigskip
Note that the action of $G$ commutes with the action of the Frobenuis morphism over the variety ${\rm{Pic}}^0_{Y_{\mathfrak p}}$, in this way, the Frobenius morphism stabilizes the eigen spaces,   $ H^1(\overline {Y}_{\mathfrak p} ,{\Bbb Q}_l)^{\chi^j_f} $.
In the following Remark we set the connection of the two following sections with \cite{D}.
\begin{rem}
 The $L$-function, $L(\chi_f,t)=
\chi_f(  \theta_{H_\mathfrak m/\Sigma_{X_\mathfrak p},T} (t  ))$, has a functional equation
$$L(\chi_f,t)=\epsilon(\chi_f,t)\cdot  L(\chi^{-1}_f,t^{-1}).$$
The term $\epsilon(\chi_f,t)$ is a "constant" function studied  in \cite{D}. By the Grothendieck theory, these $L$-functions are given by   characteristic polynomials of the Frobenius endomorphism acting on $ H^1(\overline {Y}_{\mathfrak p} ,{\Bbb Q}_l)^{\chi_f} $. In this way the above functional equation is given, (c.f.\cite[10.3.5]{D}), by
$${\rm{det}}(Id -F_\mathfrak p\cdot t)={\rm{det}}(-F_\mathfrak p\cdot t)\cdot {\rm{det}}(Id -F^\vee_\mathfrak p\cdot t^{-1}).$$
Hence,  we have  $\epsilon(\chi_f,t)=t^{2g+d-1}\cdot  P_{F_\mathfrak p}^{\chi_f}(0)$. The term  $P_{F_\mathfrak p}^{\chi_f}(0)$
  is calculated in terms of Gauss sums in \cite[Section 5]{D} and by Jacobi sums  for $g=0$ and $d=2$ in \cite{W}. In the two following sections,  by using a different point of view to that of \cite{D}, we have that $P_{F_\mathfrak p}^{\chi_f}(0)$ is given by Jacobi sums.
 \end{rem}

\section{ The constant term of $P_{F_\mathfrak p}^{\chi_f}(t)$}

In this section, we shall   calculate   $P_{F_\mathfrak p}^{\chi_f}(0)$. This means, by (*),  calculating
$$\sum_{[N]\in G_{\mathfrak m}}   L_{\mathfrak m}(N,2g+d-1 )\cdot \chi_f(\sigma_N) .$$
 We shall calculate $L_{\mathfrak m}(N,2g+d-1 )$ for each $[N]\in G_{\mathfrak m}$. Let $\kappa$ be a divisor of degree $2g-2$ associated with the canonical sheaf of $X_\mathfrak p$.
We have two cases for the cardinal
$$ \#H^0(X_{\mathfrak p}, N+(2g+d-1)\cdot D_1 -E).$$
 It is  $=q^g$ when $N+(2g+d-1)\cdot D_1 -E$ is     linearly equivalent to $\kappa$ and $=q^{(g-1)}$ in the other case.

\begin{lem}\label{no} Let $(L,\iota_\mathfrak m)$ be a level structure where ${\rm{deg}}(L)=2g+d-1$. If $L$ is a line bundle not isomorphic to ${\mathcal O}_{X_\mathfrak p}(\kappa+E)$, then $\#H^0_{\mathfrak m}((L,\iota_\mathfrak m))=q^{ (g-1)}$.
\end{lem}
\begin{proof}Let us denote ${\mathcal O}_{\mathfrak m} :=H^0( {X_\mathfrak p},{\mathcal O}_{X_\mathfrak p}/\mathfrak m)$ and let $\pi$ be the natural epimorphism ${\mathcal O}_{X_\mathfrak p}\to {\mathcal O}_{X_\mathfrak p}/\mathfrak m$. By taking global sections on the exact sequence
$$0\to L(-E)\to L\overset {\iota_\mathfrak m}\to {\mathcal O}_{X_\mathfrak p}/\mathfrak m\to 0$$
  we obtain the exact sequence of vector spaces
$$0\to H^0( {X_\mathfrak p},L(-E))\to H^0(  X_{\mathfrak p},L )\overset {H^0(\iota_\mathfrak m )}\to {\mathcal O}_{\mathfrak m} \to 0.$$
  Thus, we have an isomorphism of ${\Bbb F}_{p^h}$-vector spaces
$$ H^0(\iota_\mathfrak m ):H^0( {X_\mathfrak p},L )/H^0( {X_\mathfrak p},L(-E))\to {\mathcal O}_{\mathfrak m},$$
and we obtain
a section $s:{\mathcal O}_{X_\mathfrak p}\to L$ such that $H^0(\iota_\mathfrak m )(s)=\pi(1)$;
thus we have that $H^0_\mathfrak m(L,\iota_\mathfrak m)\neq 0$ and its cardinal is $H^0( {X_\mathfrak p},L(-E))=q^{(g-1)}$.
\end{proof}
 We now study the case when $L\simeq {\mathcal O}_{X_\mathfrak p}(\kappa+E)$.

 We denote by $\pi_{\overline x_0}$ the surjective morphism of modules
${\mathcal O}_{X_\mathfrak p}/{\mathfrak m}\to {\mathcal O}_{X_\mathfrak p}/{\mathfrak m}_{\overline x_0}$, with ${\mathfrak m}_{\overline x_0}$ the maximal ideal associated with $\overline x_0$. We denote $\iota_{\overline x_0} :=\pi_{\overline x_0}\cdot \iota_{\mathfrak m}$. Recall that $E=\overline x_0+\cdots+\overline x_d$.

We shall study when a level structure $(L,\iota_\mathfrak m)$,  with $L$ a line bundle  isomorphic to ${\mathcal O}_{X_\mathfrak p}(\kappa+E)$, has an $\mathfrak m$-section.

Bearing in mind that  if $\lambda \in k(\mathfrak p)^\times$  then $(L,\iota_\mathfrak m)$ and $(L,\lambda\cdot \iota_\mathfrak m)$ are isomorphic level structures,   we can   fix the morphism $\iota_{\overline x_0}:=\pi_{\overline x_0}\cdot \iota_\mathfrak m :L\to {\mathcal O}_{ X_\mathfrak p}/{\mathfrak m}\to {\mathcal O}_{X_\mathfrak p}/{\mathfrak m}_{\overline x_0}$:
For the level structures     $({\mathcal O}_{X_\mathfrak p}, \iota_\mathfrak m)$, we set  $    \iota_{\overline x_0}(1) =1$. Now, let us consider $\overline D$ an effective divisor with $\overline x_0\notin {\rm{supp}(\overline D )} $. Recall that by the natural inclusion $ {\mathcal O}_{X_\mathfrak p}\hookrightarrow   {\mathcal O}_{X_\mathfrak p}( \overline D )$ we have that $1\in {\mathcal O}_{X_\mathfrak p}( \overline D )$ and $1\notin {\mathfrak m}$.
We obtain an $\overline x_0$-level structures, $\iota_{\overline x_0}$, for $ {\mathcal O}_{X_\mathfrak p}(\overline D) $   by setting  $    \iota_{\overline x_0}(1) =1$.

 Because  $L$ is of degree $2g+d-1$,   we can obtain an effective divisor, $\overline D$,   with support outside $\overline x_0$  such that
$L$ is isomorphic to ${\mathcal O}_{X_\mathfrak p}(\overline D)$. We consider a $\overline D $ linearly equivalent to $\kappa+E$, with $\overline x_0\notin \text{sup}(\overline D)$, and we set $\iota_{\overline x_0}$ for $L={\mathcal O}_{X_\mathfrak p}(\overline D) $.

 We denote by  $E'$ the divisor $\overline x_1+\cdots+\overline x_d$, by $\mathfrak m'$ the ideal associated with $E'$, and ${\mathcal O}_{\mathfrak m'} :=H^0( {X_\mathfrak p},{\mathcal O}_{X_\mathfrak p}/\mathfrak m')$.

 Because  $   {\mathcal O}_{\mathfrak m }=  k(\overline x_0)\times {\mathcal O}_{\mathfrak m'}$, for each  $\mathfrak m $-level structure $(L,\iota_\mathfrak m)$  we have that $\iota_\mathfrak m=\iota_{\overline x_0}\times \iota_{\mathfrak m'}  $.

\begin{propo} We have  isomorphisms of vector spaces
$$H^0( {X_\mathfrak p},L )/H^0( {X_\mathfrak p},L(-E))\simeq H^0( {X_\mathfrak p},L )/H^0( {X_\mathfrak p},L(-E'))\simeq   {\mathcal O}_{\mathfrak m'}.$$

 Moreover, we have $H^0(\iota_{\overline x_0})(H^0( {X_\mathfrak p},L(-E)))=0$.
 \end{propo}
 \begin{proof}
 The first isomorphism is a consequence of the equality
 $$H^0( {X_\mathfrak p},L(-E))=H^0( {X_\mathfrak p},L(-E')).$$
Because of the Riemann-Roch Theorem $  H^1( {X_\mathfrak p},L(-E'))=0$, and thus  from $H^0(\iota_{\mathfrak m'})$ we obtain the second isomorphism. The last assertion follows because     $x_0\in {\rm{sup}}(E)$.
\end{proof}
We denote the composition of these two isomorphisms by $[\iota_{\mathfrak m'}]$.

 We now choose an  ${\mathfrak m'}$-level structure for  $L$, $(L,  \iota')$. We consider the ${\Bbb F}_{p^h}$-linear form
$$\omega:=H^0(\iota_{\overline x_0})\cdot [\iota']^{-1}:{\mathcal O}_{\mathfrak m'}\simeq H^0( {X_\mathfrak p},L )/H^0( {X_\mathfrak p},L(-E ))\to k(\overline x_0)={\Bbb F}_{p^h}.$$

By considering the standard basis for the ${\Bbb F}_{p^h}$-vector space ${\mathcal O}_{\mathfrak m'}=k(\overline  x_1)\times\cdots\times k(\overline  x_d)$,  we have that $\omega(z_1,\cdots,z_d)=\lambda_1\cdot  z_1+\cdots +\lambda_d\cdot  z_d$. Since
${\rm{deg}}(L)=2g+d-1$ and
$${\rm{Ker}}(\omega)=[\iota'](H^0( {X_\mathfrak p},L(-\overline x_0 )/H^0( {X_\mathfrak p},L(-E))),$$
we have that for $1\leq i \leq d,$
$$H^0( {X_\mathfrak p},L(-\overline x_0 )/H^0( {X_\mathfrak p},L(-E ))\neq H^0( {X_\mathfrak p},L(-\overline x_0-\overline x_i )/H^0( {X_\mathfrak p},L(-E )).$$
   Thus, $\lambda_1\cdots  \lambda_d\neq 0$. Accordingly, by considering   $(-\lambda_1^{-1},\cdots,-\lambda_d^{-1})\cdot  \iota'$ instead of $\iota'$, we can assume that $\lambda_1=\cdots =\lambda_d=-1$.

Note  that if $(L,\iota_{\mathfrak m}')$ is an $\mathfrak m'$-level structure, then there exists $z\in {\mathcal O}_{\mathfrak m'}^\times=k(x_1)^\times\times \cdots \times k(x_d)^\times$ with $\iota_{\mathfrak m'}= z\cdot\iota' $.

\begin{lem}\label{om} By using  the above notations, $H^0_{\mathfrak m}((L,\iota_{\overline x_0}\times z\cdot \iota'))\neq 0$ if and only if $\omega(z^{-1})=1$.
\end{lem}
\begin{proof}  If there exists a section $s$ of $L$ such that the diagram
$$ \xymatrix{ {\mathcal O}_{X_\mathfrak p}
  \ar[r]^s\ar[dr]^{\pi'}     & L\ar[d]^{z\cdot \iota'}
   \\
     &
    {\mathcal O}_{X_\mathfrak p}/\mathfrak m'}
   $$
is commutative, where $\pi'$ is   the natural epimorphism, then  we have that the class of the section $s$ in
 $  H^0( {X_\mathfrak p},L )/H^0( {X_\mathfrak p},L(-E'))$  is
$[z\cdot \iota']^{-1}(1)$.

Moreover, since the diagram
$$  \xymatrix {{\mathcal O}_{X_\mathfrak p}
  \ar[r]^{s}\ar[dr]^{\pi_{\overline x_0 }}      & L\ar[d]^{\iota_{\overline x_0}}
   \\
       & k(\overline x_0)}$$
 must also be   commutative,  we have that $H^0(\iota_{\overline x_0})(s)=1$. Thus,
 $$1=H^0(\iota_{\overline x_0})([z\cdot \iota']^{-1}(1))=H^0(\iota_{\overline x_0})([\iota']^{-1}(z^{-1}))=\omega(z^{-1}) .$$

 Reciprocally, it suffices to consider a section $s$ of $L$  within the class   $[z\cdot \iota']^{-1}(1)\in H^0( {X_\mathfrak p},L  )/H^0( {X_\mathfrak p},L(-E') )$.
 \end{proof}

 We denote by $ M $   a divisor of degree $0$ with support outside $T$ such that the $\mathfrak m$-level structure associated with  $M+(2g+d-1)\cdot D_1$ is   $({\mathcal O}_{X_\mathfrak p}(\kappa+E),\iota_{x_0}\times \iota')$, and we recall that $\sigma_M $ is the element of the Galois group of the extension $H_\mathfrak m/\Sigma_{X_\mathfrak p} $ given, via the Artin symbol, by the class
 $[M]  \in  \frac{ I^0_T}{(\Sigma^{\mathfrak m}_{X_\mathfrak p})^\times\cdot O_T^\times}$.

We  denote  $z^{-1}\cdot M$ instead of $\eta(z^{-1})\cdot [M]$.

 \begin{lem}\label{cost}
We have
   $$P_{F_\mathfrak p}^{\chi_f}(0)=\sum_{[N]\in G_{\mathfrak m}}   L_{\mathfrak m}(N,2g+d-1 )\cdot \chi_f(\sigma_N)=$$ $$=\sum_{z\in {\mathcal O}_{\mathfrak m'}^\times, \, \omega(z )=1 }q^{   g} \cdot\chi_f(\sigma_{z^{-1}\cdot M})+
  \sum_{\underset{  N+(2g+d-1)D_1\nsim \kappa+E} {[N]\in G_{\mathfrak m}}}q^{  g-1 }\cdot \chi_f(\sigma_{N}).$$
 \end{lem}
 \begin{proof} To prove the Lemma, we   bear  in mind that  the set of  classes of level structures for the line bundle
$${\mathcal O}_{X_\mathfrak p}(\kappa+E-(2g+d-1)D_1)$$
is given by
 $ \{\eta(z^{-1})\cdot [M]\}_{z\in {\mathcal O}_{\mathfrak m'}^\times } \subset \frac{ I^0_T}{(\Sigma^{\mathfrak m}_{X_\mathfrak p})^\times\cdot O_T^\times}$ and Lemmas \ref{no} and \ref{om}.
\end{proof}


\section{Jacobi sums}
In the first part of this section we follow  \cite{W}. Let us consider $z_1,\cdots,z_d\in k(\mathfrak p)^\times={\Bbb F}_{p^h}$.
We consider the Jacobi sum
$$J_{(a)}({\mathfrak p}):=(-1)^{d+1}\underset{\underset {z_1,\cdots,z_d  \quad {\rm{mod} }\, \mathfrak p}{z_1+\cdots+z_d=-1 \quad {\rm{mod }}\, \mathfrak p}}\sum \chi_{\mathfrak p}^{a_1}(z_1)\cdots \chi_{\mathfrak p}^{a_d}(z_d),$$
with $a:=(a_1,\cdots,a_d)$. The positive integers $a_1,\cdots,a_d$ are defined in the introduction. We denote
$$\theta(a):=\underset{\underset {t{\, mod \, m}} {(t,m)=1}}\sum [\sum^d_{i=1}<\frac{t \cdot a_i}{m}>]\sigma^{-1}_{-t}.$$
The map ${\mathfrak p}\to J_{(a)}({\mathfrak p})$ defines a Hecke character for the cyclotomic field $K$
and the ideal generated by $J_{(a)}({\mathfrak p})$ within ${\Bbb Z}[\epsilon_m]$ is ${\mathfrak p}^{\theta(a)}$.

Moreover, $\vert J_{(a)}({\mathfrak p})\vert^2=p^{h(s-2)}$, $s$  being the cardinal of the subset of the integers $a_1,\cdots,a_d, a_1+\cdots+a_d $, which are $\neq 0$ (mod $m$). In this case,    $s=d+1$.

 \begin{thm}\label{j}
We have
 $$P_{F_\mathfrak p}^{\chi_f}(0)=(-1)^{d+1}\chi_f(\sigma_M)\cdot q^{  g}\cdot J_{(a)}({\mathfrak p}).$$
 \end{thm}
 \begin{proof}Bearing in mind the end of  Section 2, one can consider $z:= (z_0,z_1,\cdots,z_d)\in k(\overline x_0)\times\cdots\times k(\overline x_d))^\times/k({\mathfrak p})^\times $  with $z_0=1$. Moreover, we denote
 $ \sigma_{ z^{-1}\cdot M}:= \sigma_{\eta(1,z^{-1}_1,\cdots,z^{-1}_d)}\cdot \sigma_{  M}. $
We have that
 $$\chi_f(\sigma_{z^{-1}\cdot M})=\chi^{a_1}_{\mathfrak p}(z_1)\cdots \chi^{a_d}_{\mathfrak p}(z_d)\cdot \chi_f(\sigma_{  M}). $$

Bearing in mind that  $\omega(z_1,\cdots,z_d)=-z_1-\cdots -z_d$,  we have that
$$(**)(-1)^{d+1}\chi_f(\sigma_M)\cdot q^{  g}\cdot J_{(a)}({\mathfrak p})=\sum_{z\in {\mathcal O}_{\mathfrak m'}^\times, \, \omega(z)=1 }q^g \cdot \chi_f(\sigma_{z^{-1}\cdot M})=$$ $$=\sum_{z\in {\mathcal O}_{\mathfrak m'}^\times, \, \omega(z)=1 }q^g \cdot\chi_f(\sigma_{z^{-1}\cdot M}) +\sum_{[N]\in G_{\mathfrak m}}q^{g-1}\cdot \chi_f(\sigma_{N})$$
because $\sum_{[N]\in G_{\mathfrak m}}q^{g-1}\cdot \chi_f(\sigma_{N})=0$.

We have that  for $\mu\in k(\mathfrak p)^\times$
 $$\sum_{z\in {\mathcal O}_{\mathfrak m'}^\times, \, \omega(z)=\mu }  \chi_f(\sigma_{z^{-1}    \cdot M})=$$
 $$=\chi_{\mathfrak p}^{a_1+\cdots+a_d}(\mu)\cdot
 \underset{\underset {z_1,\cdots,z_d  \quad {\rm{mod }}\, \mathfrak p}{z_1+\cdots+z_d=-1 \quad {\rm{mod }}\, \mathfrak p}}\sum \chi_{\mathfrak p}^{a_1}(z_1)\cdots \chi_{\mathfrak p}^{a_d}(z_d)\cdot \chi_f(\sigma_M),$$
 therefore the sum over all $\mu\in k(\mathfrak p)^\times$ of the above sum  is $0$, because $ a_1+\cdots+a_d\neq 0 {\rm{\, mod}} \,m$, and
  thus, together with Lemma \ref{cost} we obtain
$$P_{F_\mathfrak p}^{\chi_f}(0)=\sum_{[N]\in G_{\mathfrak m}}   L_{\mathfrak m}(N,2g+d-1 )\cdot  \chi_f(\sigma_N )=$$
$$=\sum_{z \in {\mathcal O}_{\mathfrak m'}^\times, \, \omega(z)=1 }q^g \cdot \chi_f(\sigma_{z^{-1}\cdot M})+ \sum_{\mu\in {\Bbb F}^\times_\mathfrak p} q^{(g-1)}\sum_{z\in {\mathcal O}_{\mathfrak m'}^\times, \, \omega(z)=\mu }  \chi_f(\sigma_{z^{-1}\cdot M}) +$$
$$+\sum_{\underset {N+(2g+d-1)D_1\nsim \kappa+E}{[N]\in G_{\mathfrak m}}}q^{(g-1)}\cdot \chi_f(\sigma_{N}).$$

Since  the set of  classes of level structures for the line bundle associated with $M$
 is given by
 $ \{\eta(z^{-1})\cdot [M]\}_{z\in {\mathcal O}_{\mathfrak m'}^\times } \subset G_\mathfrak m$,  the last sum is equal to
$$\sum_{z\in {\mathcal O}_{\mathfrak m'}^\times, \, \omega(z)=1 }q^g \cdot\chi_f(\sigma_{z^{-1}\cdot M}) +\sum_{[N]\in G_{\mathfrak m}}q^{(g-1)}\cdot \chi_f(\sigma_{N})$$
 and we conclude by the equality $(**)$.
 \end{proof}

\begin{rem} In this Remark we use the formula (*) written
in   section 3.

  We have
$$  {\rm{det}} (F_{\mathfrak p}) =(-1)^{2g+d-1}\sum_{[N]\in G_{\mathfrak m}}   L_{\mathfrak m}(N,2g+d-1 )\cdot \chi_f(\sigma_N).$$

Thus, by    \cite{W} and Theorem \ref{j}
$${\mathfrak p}\to  \chi_f(\sigma_M)^{-1}  {\rm{det}} (F_{\mathfrak p} )=  q^g\cdot J_{(a)}({\mathfrak p}),$$ with ${\mathfrak p}\in {\rm{Spec}} (A)$, gives a Hecke character for $K$.

We now include a formula that relates $\chi_f$, $\chi_{\mathfrak p}$,   values  over divisors ${X_\mathfrak p}$    of elements of $\Sigma_{X_\mathfrak p}$ and Jacobi sums.

We   have  the equality,
$$\sum_{[N]\in G_{\mathfrak m}}   L_{\mathfrak m}(N,2g+d-1 )\cdot \chi_f(\sigma_N)
=\sum_{D,\, {\rm{deg}} (D)=2g+d-1}    \chi_{\mathfrak p}(f(D)), $$
where the sum is over  all the   effective divisors on $X_{\mathfrak p}$  with support outside $T$. Therefore, from Theorem \ref{j} we deduce the formula
$$\sum_{D,\, {\rm {deg}}(D)}     \chi_{\mathfrak p}(f(D))=\chi_f(\sigma_M)\cdot q^g\cdot J_{(a)}({\mathfrak p}).$$
\end{rem}
We have defined $f(D)$ at the end of the section 3.

For example, let $X$ be   a proper, smooth model  for the elliptic curve $  y^2-x\cdot(x-1) \cdot (x-\lambda)$, defined over $A:={\Bbb Z}[\epsilon_m, \frac{1}{2m}]$,  $m\neq 3$. We have that
${\rm{div}}(y)=-3\cdot\infty+(0,0)+(1,0)+(\lambda, 0)$ and $a=(-3,1,1,1)$,  thus
$$\sum_{D,\, {\rm {deg}}(D)= 4}    \chi_{\mathfrak p}(y(D))= \chi_y(\sigma_M)\cdot p \cdot J_{(a)}({\mathfrak p}),$$
where the sum is over all      effective divisors $D$ on the elliptic curve, $X_\mathfrak p$, with support outside $T:=\{\infty,(0,0),(1,0),(\lambda, 0)\}$.
\section{The $p$-rank of Jacobians}

Let $ W, S$ be    proper, smooth and geometrically irreducible curves over ${\Bbb F}_q$.
  \begin{defn} The curve $  W$    has $p$-rank $0$ when,  for the morphism $[p  ]:{\rm{Pic}}^0_{ W}\to {\rm{Pic}}^0_{ W}$, $\rm{Ker}[p  ]_{{\rm{red}}}={\rm{Spec}}({\Bbb F}_q)$ or, equivalently, when
$[p]$  is purely inseparable. Here, ${\rm{Ker}}[p  ]$ is considered as a scheme and ${\rm{Ker}}[p  ]_{{\rm{red}}}$ is the reduced scheme. Note that ${\rm{Ker}}[p  ]_{{\rm{red}}}={\rm{Spec}}({\Bbb F}_q)$ if and only if ${\rm{Ker}}[p  ]({\Bbb F})=\{0\}$.
 \end{defn}

 In the following Proposition we remember a known result.  We consider $F_q$, the $q$-Frobenius morphism, as a  ${\Bbb Q}_l$-linear application over
  the ${\Bbb Q}_l$-vector space
 $H^1(\overline W ,{\Bbb Q}_l) $.
\begin{propo}\label{z} We have that ${  W}$ has $p$-rank $0$ if and only if the characteristic polynomial of
 $F_q$, $t^{2\pi}+a_1t^{2\pi-1}+\cdots+a_{2\pi-1}t+a_{2\pi}\in {\Bbb Z}[t]$, satisfies $p\vert a_1,\cdots,p\vert a_{2\pi}$. Here, $\pi$ denotes the genus of ${  W}$.
  \end{propo}
  \begin{proof}If $p\vert a_1,\cdots,p\vert a_{2\pi}$  then  by considering $F_q$ as  an endomorphism of ${\rm{Pic}}_{ W}^0$, we have that $F_q^{2\pi}=  [p  ]\cdot \Phi$ , where $\Phi$ is also an isogeny of $\rm{Pic}_{  W}^0$. Thus $[p   ]$    is purely inseparable  because $F_q$ is purely inseparable.

  Conversely, since  ${\rm{Ker}}[p  ]_{{\rm{red}}}={\rm{Spec}}({\Bbb F}_q)$ there exists $s\in {\Bbb N}$ such that
  ${  F}_q^s({\rm{Ker}}[p  ])={\rm{Spec}}({\Bbb F}_q)$. In this way,   ${  F}_q^s=[p]\cdot \Phi$ with $\Phi$ an isogeny of  ${\rm{Pic}}_{ Z}^0$.
     Thus, the characteristic polynomial of $F_q$ is $t^{2\pi}\, {\rm{mod}}\,  p $.
  \end{proof}

 Let $  W\to   S$ be a ramified abelian covering, of group $G:={\Bbb Z}/m$.
  Let $\chi$ be a non-trivial character of $G$. We  denote
 by $p_1(t)\in {\Bbb Z}[\epsilon_{m}][t]$   the characteristic polynomial of $F_q$ as    ${\Bbb Q}_l(\epsilon_{m })$-endomorphism over the ${\Bbb Q}_l(\epsilon_{m })$-vector space $  H^1(\overline W ,{\Bbb Q}_l)^\chi$, and by $p_0(t)$ the characteristic polynomial of $F_q$ over the ${\Bbb Q}_l $-vector space $ H^1(\overline S ,{\Bbb Q}_l) $.

  \begin{lem}\label{W}  If $W$ has $p$-rank $0$ then $  p_1(0)=u\cdot p $, with $u \in {\Bbb Z}[\epsilon_{m }]$.
  \end{lem}
  \begin{proof} By     Proposition \ref{z},  $F_q^s=[p ]\cdot \Phi$, where   $\Phi$ is an isogeny of the abelian variety ${\rm{Pic}}^0_W$. The characteristic polynomial of $\Phi$ as a ${\Bbb Q}_l(\epsilon_{m })$-endomorphism of $  H^1(\overline W ,{\Bbb Q}_l)^\chi$ is
  $t^{r }+b_{r -1}t^{r -1}+\cdots+b_0\in {\Bbb Z} [\epsilon_{m }]$. Note that $F_q^s$ and $[p ]$ are ${\Bbb Q}_l(\epsilon_{m })$-endomorphisms and thus $\Phi$. Therefore, the characteristic polynomial of $F_q^s$ is $t^{r }+c_{r -1}  t^{r -1}+\cdots+b_0\cdot p^{r }\in {\Bbb Z} [\epsilon_{m }]$ and we have
   $p_1(0)^s=b_0\cdot p^{r }$. We conclude because    the primary ideal decomposition of $p$ in ${\Bbb Z}[\epsilon_{m }]$ is a product of different    prime ideals.
  \end{proof}

We consider the notations of   sections 1 and 2. For $W=Y_{\mathfrak p}$, $S= X_{\mathfrak p}$ and $\chi=\chi_f$, we have $p_1(t)=P_{F_\mathfrak p}^{\chi_f}(t)$.

Let $h$ be the least positive integer with $m\vert (p^h-1)$. The ideal   $p\cdot {\Bbb Z}[\epsilon_m,\frac{1}{m\cdot n}]$  decomposes into a product of $b:={\frac{m-1}{h}}$   prime ideals  ${\mathfrak p}_1 \cdots {\mathfrak p}_b $. The action of $({\Bbb Z}/m)^\times$ on  $\{{\mathfrak p}_1, \cdots, {\mathfrak p}_b\} $ has as isotropy $C_h$,  the unique $h$-cyclic subgroup of $({\Bbb Z}/m)^\times$.

By choosing the representative $i$ within the class $[i]$, we obtain the identification, as sets, $\{ 1,\cdots, m-1\}=({\Bbb Z}/m)^\times$. We denote by $c_h\subset \{ 1,\cdots, m-1\}$ the subset associated, by the above identification,    with $C_h$.
We consider the group quotient $({\Bbb Z}/m)^\times/C_h$, which acts transitively and without isotropy on $\{{\mathfrak p}_1, \cdots, {\mathfrak p}_b\} $. Let $e_h:=\{i_1,\cdots,i_b\}\subseteq  \{ 1,\cdots,m-1\}$ be representatives of the classes of $({\Bbb Z}/m)^\times/C_h$.

We use $<\frac{a}{m}>$ to denote the fractional part of $\frac{a}{m}$ and
$d_u:=[<\frac{u\cdot a_1}{m}>+\cdots +<\frac{u\cdot a_d}{m}>]$ with $u\in \{ 1,\cdots, m-1\}$.
The values $d_u$ are the dimensions of the  eigenspaces of $H^1(Y_\mathfrak p, {\mathcal O}_{Y_\mathfrak p} )$ for the different characters of $G$, c.f.\cite{Bo}.
We set $O_t:=\underset{u\in c_h} \sum d_{t\cdot u}$, which is  the dimension of the orbit for the action of $C_h$ on the  eigenspaces  $ H^1( {Y}_{\mathfrak p} ,{\mathcal O}_{Y_\mathfrak p} )^{\chi^j_f} $.

\begin{thm}\label{cr} 1) If  there exists $t\in e_h$ with
$O_1\neq O_t$
then   ${\rm{Pic}}^0_{Y_{\mathfrak p}}$ is not supersingular.

2) If $X_{\mathfrak p}={\Bbb P}^1$  and there exists $t\in e_h$ with
$  O_t= 0$,
 then $Y_{\mathfrak p}$ does not have   $p$-rank $0$.

\end{thm}
\begin{proof}1) By Theorem \ref{j},  the ideal within ${\Bbb Z}[\epsilon_m,\frac{1}{m\cdot n}]$ generated by the term constant of $P_{F_\mathfrak p}^{\chi_f}(t)$ is  $(q^g\cdot J_{(a)}({\mathfrak p}))$. Thus, by using
 \cite{W} (8), (9) this ideal is
$(q^g\cdot \mathfrak p^{\theta(a)})$
with
$\theta(a):=\sum_{ 1\leq u<m  }d_u\cdot \sigma^{-1}_{-u}. $

    Because
$\underset{u\in h_r} \sum d_u\neq \underset{u\in h_r} \sum d_{u\cdot t}$
for some $t\in e_r$, we have that the ideal primary decomposition of $(q^g\cdot J_{(a)}({\mathfrak p}))$ is ${\mathfrak p}_1^{n_1}\cdots {\mathfrak p}_b^{n_b}$, with $n_1\neq n_v$, for some $v$   and  we conclude. Bear in mind that ${\rm{Pic}}_{Y_\mathfrak p}^0$ is supersingular if and only if the eigenvalues of the  Frobenius morphism are $\zeta\cdot q^{1/2}$, $\zeta$ being a root of the unity.

  2) Since $g=0$, the  ideal     $P_{F_\mathfrak p}^{\chi_f}(0)\cdot {\Bbb Z}[\epsilon_m,\frac{1}{m\cdot n}]$   is  $(  J_{(a)}({\mathfrak p}))$.
Because  $  \underset{u\in h_r} \sum d_{u\cdot t}= 0$,   the ideal primary decomposition of this ideal is ${\mathfrak p}_1^{n_1}\cdots {\mathfrak p}_b^{n_b}$ with  $  n_v=0$  for some $v$. We conclude because      if $Y_{\mathfrak p}$ has   $p$-rank $0$, then $  P_{F_\mathfrak p}^{\chi_f}(0)=u\cdot p $ (c.f Lemma \ref{W}).
 \end{proof}

 \begin{lem}\label{Q} We have that $Y_{\mathfrak p}$ has $p$-rank $0$ if and only if  $X_{\mathfrak p}$ has $p$-rank $0$ and $P_{F_\mathfrak p}^{\chi^j_f}(t)=t^{r }+p\cdot Q_j(t)$, where $Q_j(t)\in {\Bbb Z}[\epsilon_m][t]$ and $deg(Q_j (t))<r $, for each $j$, $1\leq j\leq m-1$.
 \end{lem}
 \begin{proof} This follows  from   Proposition \ref{z},
bearing in mind that   the characteristic polynomial of the $p^h$-Frobenius morphism,  considered as a  ${\Bbb Q}_l$-linear application of  $H^1(\overline Y_{\mathfrak p},{\Bbb Q}_l) $, is   $p_0(t)\cdot\prod_{(j,m)=1}   P_{F_\mathfrak p}^{\chi^j_f}(t)$ and that the primary ideal decomposition of $p$ in ${\Bbb Z}[\epsilon_m,\frac{1}{m\cdot n}]$  is a product of different  prime ideals.
 \end{proof}

 In the following theorem we assume the notation of   Theorem \ref{cr} and that  the sums are over all     effective divisors $D$ on $X_\mathfrak p$  with support outside $T$.
\begin{thm}\label{in} 1) Let $p$ be a integer prime  inert  in ${\Bbb Z}[\epsilon_m, \frac{1}{m\cdot n}]$. We have that $Y_\mathfrak p$ has $p$-rank $0$ if and only if $X_\mathfrak p$ has $p$-rank $0$ and
$$\sum_{D, {\rm{deg}}(D)=l}    \chi_{\mathfrak p}(  f(D))=0\,{\rm{mod}}\,  \mathfrak p,$$
for each $l$, $1\leq l\leq 2g+d-2$. Note that $k(\mathfrak p)={\Bbb F}_{p^{m-1}}$.

2)We now assume   $g=0$ and that
$  \underset{u\in c_h} \sum d_{u\cdot t}\neq 0$
for each   $t\in e_h$. We have that $Y_\mathfrak p$ has $p$-rank $0$ if and only if
$$\sum_{D, {\rm{deg}}(D)=l}    \chi_{\mathfrak p}(  f(D))^{j}=0\,{\rm{mod}}\, \mathfrak p $$

for each   $j\in e_h$ and each $l$, $1\leq l\leq   d-2$.
 Note that in this case $k(\mathfrak p)={\Bbb F}_{p^h }$.

\end{thm}
\begin{proof} 1) From   (*) section 3, we have
$$P_{F_\mathfrak p}^{{\chi_f} }(t)=\underset{i=0}{\overset
{ 2g+d-1} \sum}\sum_{D, {\rm{deg}}(D)=i}    \chi_{\mathfrak p}( f(D))t^{2g+d-1-i}.$$
We conclude by using the   Lemma \ref{Q}
and  the fact  that by Theorem \ref{j},   $$P_{F_\mathfrak p}^{\chi_f}(0)\cdot {\Bbb Z}[\epsilon_m,\frac{1}{m\cdot n}]=(q^g\cdot J_{(a)}({\mathfrak p})).$$
Note  that in the case of $p$ being inert  in ${\Bbb Z}[\epsilon_m,\frac{1}{m\cdot n}]$ , ${\mathfrak p}= p\cdot {\Bbb Z}[\epsilon_m,\frac{1}{m\cdot n}]$  and
 $q^gJ_{(a)}({\mathfrak p})=\epsilon_m^j\cdot p^u$  for some $1\leq u,\, j  \in {\Bbb N}$. Moreover, $P_{F_\mathfrak p}^{\chi^j_f}(t)= \sigma^j(P_{F_{\mathfrak p}}^{\chi_f}(t))= 0{\rm{\, mod}}\,  {\mathfrak p}$, for each $j$ with  $1\leq j \leq  m-1$, if and only if
$P_{F_{\mathfrak p}}^{\chi_f}(t) = 0{\rm{\, mod}}\,  {\mathfrak p}$.

2) One proceeds in the same way as 1). Bearing in mind that   $g=0$ and that
 $  \underset{u\in c_h} \sum d_{u\cdot t}\neq 0$,
for each   $t\in e_h$, we have that,  $(  J_{(a)}({\mathfrak p} ))=(  \mathfrak p^{\theta(a)})\subseteq p\cdot  {\Bbb Z}[\epsilon_m,\frac{1}{m\cdot n}]$.
Note that  $P_{F_\mathfrak p}^{\chi^j_f}(t)=\sigma_{j}( P_{F_\mathfrak p}^{\chi_f}(t))$, $p\cdot {\Bbb Z}[\epsilon_m,\frac{1}{m\cdot n}]  ={\mathfrak p}_1\cdots {\mathfrak p}_b$,   and  $\{\sigma_j\}_{1\leq j \leq m-1}$ operates  transitively on ${\mathfrak p}_1,\cdots, {\mathfrak p}_b$.
 \end{proof}

 In the next Corollary, by using  part 1) of this Theorem  we obtain counterparts to the Deuring polynomial for genus greater than $1$. In \cite[5]{Bo}, it is  obtained, explicitly, the generalized Hasse-Witt matrix for $Y_{\mathfrak p}$ and  the author thus   obtained equations to study the $p$-rank of $Y_{\mathfrak p}$.

   We consider the curves defined over ${\Bbb Z}[\epsilon_m,\frac{1}{m\cdot n}]$ , with  $n$   the product of integer primes $p$ with $\#k(\mathfrak p)=p^r<d$,   $X={\Bbb P}^1 $  and  $Y$ the curve associated with the plane curve $$y^m={ x^{a_1}(x-1)^{a_2}(x-\alpha_1)^{a_3}\cdots(x-\alpha_{d-2})^{a_d}} ,$$
where the $a_1,\cdots, a_d,a_1+\cdots+a_d$ are  integers $\neq 0\,{\rm{mod}}\, m$, with $1\leq a_j<m$    and $\alpha_1,\cdots,\alpha_{d-2}$ are different elements of $k(\mathfrak p)\setminus \{0,1\}$ and $f=x^{a_1}(x-1)^{a_2}(x-\alpha_1)^{a_3}\cdots(x-\alpha_{d-2})^{a_d}$.

\begin{cor}   Let $Y_\mathfrak p$ be the proper, smooth model of   the reduction at $\mathfrak p$ of $Y$, with $p$  inert in ${\Bbb Z}[ \epsilon_m,\frac{1}{m\cdot n}]$.  Accordingly, $Y_\mathfrak p$ has $p$-rank $0$   if and only if
$$\sum_{ \underset{{\rm{deg}}(q(x))=l}  {q(x)}}[f({\rm{div}}( q(x) ))]^{(p^{m-1}-1)/m}=0\,$$
or equivalently,
$$\sum_{ \underset{{\rm{deg}}(q(x))=l}  {q(x)}}[q(0)^{a_1}q(1)^{a_2}q(\alpha_1)^{a_3}\cdots q(\alpha_{d-2})^{a_d}]^{(p^{m-1}-1)/m}=0\,$$
for each $l$, with $1\leq l \leq d-2$. Here, the sums are over all    monic polynomials  $q(x)\in {\Bbb F}_{p^{m-1}}[x]$.
\end{cor}
\begin{proof}
 It suffices to consider in   Theorem \ref{in} that  $X_\mathfrak p={\Bbb P}^1$, $T:=(x\cdot (x-1)\cdot (x-\alpha_1)\cdots (x-\alpha_{d-2}))_0\cup \{\infty\}$  and also to consider  that the effective divisors of degree $l$ on ${\rm{Spec} }({\Bbb F}_{p^{m-1}}[x])$ are given by the zero locus of monic polynomials $q(x)$ of degree $l$. To conclude it suffices bear in mind that if ${\rm{div}}( q(x) )=\sum_in_i \cdot y_i-l\cdot \infty$, then $f({\rm{div}}( q(x) ))=\prod_if(y_i)^{n_i}\cdot f(\infty)^{-l}$. By the Weil's reciprocity law, we have
$$f(\sum_in_i \cdot y_i)=(-1)^{l\cdot(a_0+\cdots+a_d)}q(0)^{a_1}q(1)^{a_2}q(\alpha_1)^{a_3}\cdots q(\alpha_{d-2})^{a_d}.$$
 Note that  $q(x)$ and $f(x)$ are monic polynomials and thus $ q^{\rm{deg}(q(x))}f^{-l}(\infty)=1$. Moreover, it is not necessary to impose to  $q(x)$  that $q(0)\cdot q(1)\cdot q(\alpha_1)\cdots q(\alpha_{d-2})\neq 0 $ because   in  the above sum   these terms are $0$.

\end{proof}
We set $\pi$ the genus of $Y_{\mathfrak p}$. By    \cite{A} 4.1.1, we have that the zeta function for $Y_{\mathfrak p}$ is:

 $$    H\cdot (\frac{   T^{2\pi }}{(1-T)(1-qT)})+
  \sum^{ \pi -1}_{i=0}     A_i\cdot T^{i}+   q\cdot  \sum^{ 2\pi -2}_{i=\pi }    B_i\cdot T^{i}, $$
with $A_i,B_i,H\in {\Bbb Z}$.
Thus, the characteristic polynomial ${\rm{det}}(t-F_{\mathfrak p})\in {\Bbb Z}[t]$  is
$   \sum^{2 \pi } _{i=\pi}     a_i\cdot t^i+    q\cdot \sum^{ \pi -1}_{i=0}    b_i\cdot t^{i}$, with $a_i, b_i \in  {\Bbb Z}$. Therefore, in the above corollary the last $\pi-1$ equations are $0$ ${\rm{mod}}\, p$.

Now, the   system of equations   of this  corollary   is a system on  ${\Bbb F}_{p^{m-1}}$ with $d-2$ variables and $d-2-(g-1)$ equations. In \cite{GlP} it is shown that the moduli space of hyperelliptic curves with   $p$-rank $0$ is $\pi-1$.

The hyperelliptic curve defined over ${\Bbb F}_{p }$
$$y^2-x(x-1)(x-\alpha_1)(x-\alpha_2)(x-\alpha_3)$$
 has $p$-rank $0$  with $p\geq 5$, if and only if
 $$(1)\sum_{a \in {\Bbb F}_p }[a(1+a )\prod_{i=1}^3( \alpha_i+a)]^{(p-1)/2}=0$$
$$(2)\sum_{a,b\in {\Bbb F}_p }[a(1+b+a )\prod_{i=1}^3(\alpha_i^2+ \alpha_i\cdot b+a)]^{(p-1)/2}=0$$

Let $A$ be   the Cartier-Manin matrix, which is the $2\times 2$-matrix whose entry $a_{ij}$ is the coefficient $c_{ip-j}$ of $x^{ip-j}$ in the polynomial $$[x(x-1)(x-\alpha_1)(x-\alpha_2)(x-\alpha_3)]^{(p-1)/2} .$$
In \cite[Theorem 3.1]{Y} the conditions   for this hyperelliptic curve to be of $p$-rank $0$  are given by the four algebraic equations obtained from   $A^2\ =0$. Therefore, this hyperelliptic curve has $p$-rank $0$ if and only if ${\rm{Trace}}(A)=c_{p-1}+c_{2p-2}=0$ and ${\rm{det}}(A)=c_{p-1} c_{2p-2}-c_{p-2} c_{2p-1}=0$.

We  are now going to prove that these equations match  the  equations below  (1) and (2). Equation (1) is equal to $\underset{\alpha \in {\Bbb F}_p  }\sum f(\alpha)^{\frac{p-1}{2}}$. Bearing in mind that $\underset{\alpha \in {\Bbb F}_p  }\sum \alpha^s=0$ if and only if $s$ is a multiple of $p-1$, we have that $(1)=-c_{p-1}-c_{2p-2}$.

Equation (2) is equal to $$1/2 (\underset{\alpha \in {\Bbb F}_{p^2}  }\sum  [f(\alpha)f(\alpha^p)]^{\frac{p-1}{2}}+
\underset{(\alpha,\beta) \in {\Bbb F}_{p }\times  {\Bbb F}_{p } }\sum [f(\alpha)f(\beta)]^{\frac{p-1}{2}}).$$
Now,   the only terms that we must consider in the first summand are the coefficients of $x^{p^2-1}$ and $x^{2(p^2-1)}$ in $[f(x)f(x^p)]^{\frac{p-1}{2}}$.  We consider   $x^ix^{jp}$, $0<i,j<5{\frac{p-1}{2}}$ with $i+jp$ equal to either $p^2-1   $ or $2(p^2-1)$. We have the following posibilities: $i= p-2$  $ j=2 p-1 $
 $i=2 p-1  $  $ j=p-2 $, $i=j=p-1$ and $i=j=2p-2$. Thus, $\underset{\alpha \in {\Bbb F}_{p^2}  }\sum  [f(\alpha)f(\alpha^p)]^{\frac{p-1}{2}}=  -2c_{2 p-1}c_{  p-2}-c_{p-1}^2 -c_{2p-2}^2 $.

For the  sum $\underset{(\alpha,\beta) \in {\Bbb F}_{p }\times  {\Bbb F}_{p } }\sum [f(\alpha)f(\beta)]^{\frac{p-1}{2}} $ we must consider in  $[f(x)f(y)]^{\frac{p-1}{2}}$ the coefficients of $x^i y^j$ with $i,j\in\{p-1,2p-2\}$. In this way, this
 sum is equal to $ 2  c_{2 p-2}c_{  p-1}+c_{p-1}^2 +c_{2p-2}^2 $. Thus  $ (2)=  -c_{2 p-1}c_{  p-2}+ c_{2 p-2}c_{  p-1}$ and we conclude.

The next Corollary is an application of part 2) of     Theorem \ref{in} when $r=1$. Let $Y$ be the above plane curve.
\begin{cor}   Let us consider
$$ [<\frac{t\cdot  a_1}{m}>+\cdots +<\frac{t\cdot a_d}{m}>]\neq 0$$
for each $t$, $1\leq t \leq m-1$. We have that    $Y_\mathfrak p$ has $p$-rank $0$ if and only if
$$\sum_{ \underset{{\rm{deg}}(q(x))=l}  {q(x)}}[f({\rm{div}}( q(x) ))]^{j(p^{m-1}-1)/m}=0\,$$
or equivalently,

$$\sum_{ \underset{{\rm{deg}}(q(x))=l}  {q(x)}}[q(0)^{a_1} q(1)^{ a_2}q(\alpha_1)^{ a_3}\cdots q(\alpha_{d-2})^{a_d}]^{j(p -1)/m}=0\,{\rm{mod}}\,  {\mathfrak p},  $$
for each $l$ and $j$ with $1\leq l \leq d-2$ and  $1\leq j \leq m-1$, respectively. The sums are over all    monic polynomials  $q(x)\in {\Bbb F}_{p }[x]$.

\end{cor}
\begin{proof}
 This is an application   of part 2) of    Theorem \ref{in}. With the values $f(D)$ one proceeds in the same way as in the above Corollary  by using Weil's reciprocity law.

\end{proof}

As  an example of these corollaries    we take  $Y\equiv y^3-x\cdot(x-1)\cdot (x-\alpha)^2=0$.

Let $p$ be an integer prime $\geq 5$. We have that $p$ in ${\Bbb Z}[\epsilon_3,\frac{1}{6}]$ either splits completely when $p\equiv 1 \text{ mod } 3 $  or is inert when $p\equiv -1 \text{ mod } 3 $.

We have that for all $t\in \{1,2\}$,  $ [<\frac{t }{3}>+<\frac{t }{3}>+ <\frac{2t }{3}>]\neq 0$. Hence,  the desingularization of this curve has $p$-rank $0$ over ${\Bbb F}_p$, with $p\equiv 1 \text{ mod } 3 $, if and only if
$$\sum_{a\in {\Bbb F}_p} [a\cdot(a-1)\cdot (a-\alpha)^2]^{(p-1)/3}=0,\quad \sum_{a\in {\Bbb F}_p}[a^2\cdot(a-1)^2\cdot (a-\alpha)^4]^{(p-1)/3}=0.$$

For $p\equiv -1 \text{ mod } 3 $ this curve has $p$-rank $0$ over ${\Bbb F}_{p^2}$ if and only if
$$\sum_{a\in {\Bbb F}_{p^2}} [a\cdot(a-1)\cdot (a-\alpha)^2]^{({p^2}-1)/3}
=  0.$$

As an example of this corollary one proves  that the desingularization of the above curve, defined over ${\Bbb F}_7$, has $7$-rank $\neq 0$. Moreover, by using the upper bounds of \cite[Theorem 1.1]{E}, this curve has $p$-rank $\leq 2$.

{\bf{Acknowlegment}}: I would like to express my gratitude to Prof. Jes\'us Mu\~noz D\'iaz. I would like to thank the referees  for noting a mistake in Lemma 6.3,   the important reference \cite{D},   for noting the adequate references to this article and the detailed and hard work made   to substantially improve the original article.

\vskip1.5truecm { \'Alvarez V\'azquez, Arturo}\newline {\it
e-mail: } aalvarez@usal.es


\vskip2truecm
\begin{thebibliography}{99}


\bibitem[A]{A} Anderson, G.  ''Rank one elliptic modules $A$-modules and $A$-harmonic series'',
 Duke Mathematical Journal. {\bf 73} (1994), pp.491--542


\bibitem[Be]{Be} Beckmann, S.  ''Ramified primes in the field of moduli of branched coverings of curves''
J. Algebra {\bf 125}  (1989), no. 1, pp.236-255.



\bibitem[Bo]{Bo} Bouw,  I.  ''The p-rank of ramified covers of curves'',
Compositio Math. {\bf 126}  (2001), no. 3, pp.295-322.





\bibitem[D]{D} Deligne, P.  ''Les constantes des \'equations fonctionnelles des fonctions L'', Proc.
Antwerpen Conference, vol. 2; Lecture Notes in Math. {\bf 349} (Springer-Verlag
1973), pp.501--597

\bibitem[E]{E} Elkin, A.  ''The rank of the Cartier operator on cyclic covers of the projective line''
J. Algebra {\bf 327}  (2011), pp.1-12.



\bibitem[G]{G} Gonz\'alez, J.  ''Hasse-Witt matrices for the Fermat curves of prime degree'',
Tohoku Math. J. (2) {\bf 49}  (1997), no. 2, pp.149-163.

\bibitem[GlP]{GlP} Glass, D; Pries, R.  ''Hyperelliptic curves with prescribed p-torsion'',
Manuscripta Math. {\bf 117}  (2005), no. 3, pp.299-317.



 \bibitem[N]{N} Noboru, A.  ''On supersingular cyclic quotients of Fermat curves'',
Comment. Math. Univ. St. Pauli . {\bf 57} (2008),  no. 1, pp.65-90





\bibitem[S]{S} Serre, J.P.  ''Groupes alg\'ebriques et corps de classes'', Hermann
 (1959)
\bibitem[T]{T} Tate, J.  ''Les Conjectures de Stark sur Les fonction L d'Artin en s=0'',  Birkhauser, Boston, 1984.



\bibitem[W]{W} Weil, A.  ''Jacobi sums as "Grossencharaktere" '',
Trans. Amer. Math. Soc. {\bf 73} (1952),  pp. 487-495

\bibitem[Y]{Y} Yui, N.  ''On the Jacobian varieties of hyperelliptic curves over fields of characteristic $p>2$.'',
 J. Algebra 52 (1978), no. 2, pp.378-410.









\end{thebibliography}
\end{document}